\documentclass[]{interact}

\usepackage{epstopdf}% To incorporate .eps illustrations using PDFLaTeX, etc.
\usepackage[caption=false]{subfig}% Support for small, `sub' figures and tables

\usepackage[numbers,sort&compress]{natbib}% Citation support using natbib.sty
\bibpunct[, ]{[}{]}{,}{n}{,}{,}% Citation support using natbib.sty
\makeatletter% @ becomes a letter
\def\NAT@def@citea{\def\@citea{\NAT@separator}}% Suppress spaces between citations using natbib.sty
\makeatother% @ becomes a symbol again
\usepackage{anyfontsize}
\theoremstyle{plain}% Theorem-like structures provided by amsthm.sty
\newtheorem{theorem}{Theorem}[section]
\newtheorem{lemma}[theorem]{Lemma}
\newtheorem{corollary}[theorem]{Corollary}
\newtheorem{proposition}[theorem]{Proposition}

\theoremstyle{definition}
\newtheorem{definition}[theorem]{Definition}
\newtheorem{remark}[theorem]{Remark}
\newtheorem{example}[theorem]{Example}

\usepackage{xcolor}

\theoremstyle{remark}

\begin{document}
	
	%\articletype{ARTICLE TEMPLATE}% Specify the article type or omit as appropriate

	\title{The Moore-Penrose inverses of unbounded closable operators and the direct sum of closed operators in Hilbert spaces}
  \author{\name {Arup Majumdar \thanks{Arup Majumdar (corresponding author). Email address: arupmajumdar93@gmail.com}and P. Sam Johnson \thanks{ P. Sam Johnson. Email address: sam@nitk.edu.in}} \affil{Department of Mathematical and Computational Sciences, \\
			National Institute of Technology Karnataka, Surathkal, Mangaluru 575025, India.}}
  
	\maketitle

	\begin{abstract}
		In this paper, we present some interesting results to characterize the Moore-Penrose inverses of unbounded closable operators and the direct sum of closed operators in Hilbert spaces.
	\end{abstract}
	
	\begin{keywords}
		Moore-Penrose inverse, closable operator, closed operator.
	\end{keywords}
      \begin{amscode}47A05; 47B02.\end{amscode}
    \section{Introduction}
    The Moore-Penrose inverse is a fundamental concept in linear algebra and functional analysis, providing a generalized notion of inverses for matrices and linear operators that may not have traditional inverses. The Moore-Penrose inverse is named after E. H. Moore and Roger Penrose, who independently introduced the concept in 1920 and 1955, respectively. Moreover, Yu. Ya. Tseng defines the Moore-Penrose (Generalized) inverses of densely defined linear operators on Hilbert spaces in \cite{MR0029479, MR0031192, MR0031191}. The Moore-Penrose inverses of closed operators are extensively studied in \cite{MR0136996, MR0473881}. In \cite{MR0473881}, a series representation for the Moore-Penrose inverse of a closed linear operator has been established. In paper \cite{majumdar2024hyers}, authors have introduced the Moore-Penrose inverses of closable operators with decomposable domains while establishing Lemma 2.10 \cite{majumdar2024hyers}. The Moore-Penrose inverse has been extensively investigated over the years due to its usefulness in Optimization Problems, Singular Value Decomposition, Signal Processing and Control Theory. This paper delves into the exploration of the Moore-Penrose inverses of unbounded closable operators in Hilbert spaces, elucidating an example of a closable operator which is not closed in Section 2. Section 3 investigates some properties of the Moore-Penrose inverses of the direct sum of closed operators in Hilbert spaces.

From now on, we shall mean $H$, $K$, $H_{i}$, $K_{i}$ ($i = 1,2,\dots, n$) as Hilbert spaces. The specification of a domain is an essential part of the definition of an unbounded operator, usually defined on a subspace.  Consequently, for an operator $T$, the specification of the subspace $D$ on which $T$ is defined, called the domain of $T$, denoted by $D(T)$, is to be given. The null space and range space of $T$ are denoted by $N(T)$ and $R(T)$, respectively. $W_{1}^{\perp}$ denotes the orthogonal complement of a set $W_{1}$ whereas $W_{1} \oplus^{\perp} W_{2}$ denotes the orthogonal direct sum of the subspaces $W_{1}$ and $W_{2}$ of a Hilbert space. Moreover,  $T\vert_{W}$ denotes the restriction of $T$ to a subspace $W$ of a specified Hilbert space.  We call $D(T)\cap N(T)^\perp $, the carrier of $T$ and it is denoted by $ C(T)$. $T^{*}$ denotes the adjoint of $T$, when $D(T)$ is densely defined in the specified Hilbert space. Here, $P_{V}$ is the orthogonal projection on the closed subspace $V$ in the specified Hilbert space and the set of bounded operators from $H$ into $K$ is denoted by $B(H, K)$.
	For the sake of completeness of exposition, we first begin with the definition of closed and closable operators.

	\begin{definition}
		Let $ T$ be an operator from a Hilbert space $H$ with domain $D(T) $ to a Hilbert space $K$. If the graph of $T$ defined by 
		$$ G(T)=\left\{(x,Tx): x\in D(T)\right\} $$ is closed in $H\times K $, then $T$ is called a closed  operator. Equivalently, $T $ is a closed operator if $ \{x_n \}$ in $D(T) $ such that $ x_n\rightarrow x $ and $ Tx_n\rightarrow y$ for some $ x\in  H,y\in   K $,  then $ x\in  D(T) $ and $ Tx=y $. That is, $G(T)$ is  a closed subspace of $H\times K$ with respect to the graph norm $\|(x, y)\|_T=(\|x\|^2+\|y\|^2)^{1/2}$. It is easy to show that the graph norm  $\|(x, y)\|_T$ is equivalent to the norm $\|x\|+\|y\|$. 	 We note that,  for any densely defined closed operator $T$,  the closure of $C(T)$, that is,  $ \overline{C(T)}$ is $ N(T)^\perp$. 
		We say that $S$ is an extension of $ T $ (denoted by $ T\subset S$) if $ D(T)\subset D(S) $ and $ Sx=Tx $ for all $ x\in D(T)$.
		
		An operator $ T $ is said to be closable 	
		if $T$ has a closed extension. It follows that $T$ is closable if the closure $\overline {G(T)}$ of $G(T)$ is a graph of an operator. 	It is also possible for a closable operator to have many closed extensions.  Its minimal closed extension is denoted by $ \overline{T}$. That is, every closed extension of $ T $ is also an extension of $ \overline{T}$.
		
	\end{definition}

	\begin{definition}
		Let $T$ be a closed operator from  $D(T) \subset H$ to  $K$. The generalized inverse of $T$ is the map $T^{\dagger}: R(T) \oplus^{\perp} R(T)^{\perp} \to H$ defined by
		\begin{equation}\label{equ 1}
			T^{\dagger} y = 
			\begin{cases}
				({T}\vert_{C(T)})^{-1}y  & \text{if} ~  y\in R(T)\\
				0    & \text{if}  ~ y\in R(T)^{\perp}.
			\end{cases}
		\end{equation}
    \end{definition}
\noindent  It can be shown that $T^{\dagger}$ is closed when $T$ is closed. 
	\begin{definition}\cite{majumdar2024hyers}
		A linear operator $T$ from $D(T) \subset H$ to  $K$ has a decomposable domain if $D(T) = N(T) \oplus^{\perp}  C(T)$. The generalized inverse of $T$ is the map $T^{\dagger}: R(T) \oplus^{\perp} R(T)^{\perp} \to H$ defined as follows:
		\begin{equation}\label{equ 2}
			T^{\dagger} y = 
			\begin{cases}
				(T\vert_{C(T)})^{-1}y  & \text{if} ~  y\in R(T)\\
				0    & \text{if}  ~ y\in R(T)^{\perp}.
			\end{cases}
		\end{equation}
  \begin{theorem}\label{thm 1.4}\cite{MR0396607}
  Let $T$ be a densely defined closed operator from $D(T) \subset H$ into $K$. Then the following statements hold good:
  \begin{enumerate}
  \item $T^{\dagger}$ is a closed operator from $K$ into $H$;
  \item $D(T^{\dagger}) = R(T) \oplus^{\perp} N(T^{*})$; $N(T^{\dagger}) = N(T^{*})$;
  \item $R(T^{\dagger}) = C(T)$;
  \item $T^{\dagger}Tx = P_{\overline{R(T^{\dagger})}}x, \text{ for all } x\in D(T)$;
  \item $TT^{\dagger}y = P_{\overline{R(T)}}y, \text{ for all } y\in D(T^{\dagger})$;
  \item $(T^{\dagger})^{\dagger} = T$;
  \item $(T^{*})^{\dagger} = (T^{\dagger})^{*}$;
  \item $N((T^{*})^{\dagger})= N(T)$;
  \item $(T^{*}T)^{\dagger} = T^{\dagger}(T^{*})^{\dagger}$;
  \item $(TT^{*})^{\dagger} = (T^{*})^{\dagger} T^{\dagger}$.
  \end{enumerate}
  \end{theorem}
	\end{definition}
 \section{Characterization of the Moore-Penrose inverses of unbounded closable operators:}
 Throughout this section, we assume that $T$ is a closable operator from $D(T) \subset H$ into $K$ with the decomposable domain $D(T) = N(T) \oplus^{\perp} C(T)$ and $N(\overline{T}) = \overline{N(T)}$.
 \begin{proposition}\label{pro 2.1}
 Let $T$ be a closable operator from $ D(T) \subset H$ into $K$. Then $T^{\dagger}$ is also closable and $\overline{T^{\dagger}} \subset (\overline{T})^{\dagger}$.
 \end{proposition}
 \begin{proof}
 $(\overline{T})^{\dagger}$ is closed with domain $D((\overline{T})^{\dagger}) = R(\overline{T}) \oplus^{\perp} R(T)^{\perp} \supset D(T^{\dagger})$. Moreover, $C(T) \subset C(\overline{T})$ because of $N(\overline{T}) = \overline{N(T)}$. Thus,
 \begin{align*}
 T^{\dagger}(y) &= 0 = (\overline{T})^{\dagger}(y), \text{ when } y \in R(T)^{\perp}\\
&= (T\vert_{C(T)})^{-1}(y) = (\overline{T}\vert_{C(\overline{T})})^{-1}(y) = (\overline{T})^{\dagger}(y), \text{ when } y \in R(T).
 \end{align*}
 It confirms that $T^{\dagger}$ is closable and $\overline{T^{\dagger}} \subset (\overline{T})^{\dagger}$.
 \end{proof}
 \begin{theorem}\label{thm 2.2}
 Let $T$ be a densely defined closable operator from $D(T) \subset H$ into $K$. Then $(T^{\dagger})^{*} = (T^{*})^{\dagger}$.
 \end{theorem}
 \begin{proof}
 $(T^{\dagger})^{*}$ exists because of the denseness of domain $D(T^{\dagger})$. The denseness of $D(T)$ confirms that the existence of $T^{*}$. From Proposition \ref{pro 2.1} and the closeness of $\overline{T}$ guarantee that the relation $(T^{*})^{\dagger} \subset (T^{\dagger})^{*}$. Now we will show that $D((T^{\dagger})^{*} )\subset D((T^{*})^{\dagger})$. Let us consider $z \in D((T^{\dagger})^{*})$. Then $f(y) = \langle T^{\dagger}y, z \rangle, \text{ for all }y \in D(T^{\dagger}),$ is continuous. The denseness of $D(T^{\dagger})$ and $D(T^{\dagger}) \subset D((\overline{T})^{\dagger})$ say that $g(y)= \langle \overline{T}^{\dagger}y, z\rangle, \text{ for all } y\in D((\overline{T})^{\dagger}),$ is also continuous. So, $z \in D(T^{*})^{\dagger}$ which confirms that $ D((T^{\dagger})^{*}) \subset D(T^{*})^{\dagger}$. Therefore, $(T^{\dagger})^{*} = (T^{*})^{\dagger}$.
 \end{proof}
 
 \begin{theorem}\label{thm 2.3}
 Let $T$ be a densely defined closable operator from $D(T) \subset H$ into $K$. Then $T \subset (T^{\dagger})^{\dagger} \subset \overline{T}$.
 \end{theorem}
 \begin{proof}
 $T^{\dagger}$ has decomposable domain because
 \begin{align*}
  D(T^{\dagger}) &= R(T) \oplus^{\perp} R(T)^{\perp}\\
   &= \overline{R(T)} \cap D(T^{\dagger}) \oplus^{\perp} N(T^{\dagger})\\
   &= N(T^{\dagger})^{\perp} \cap D(T^{\dagger}) \oplus^{\perp} N(T^{\dagger}) \\
   &= C(T^{\dagger})\oplus^{\perp} N(T^{\dagger}).
 \end{align*} 
 Thus, $(T^{\dagger})^{\dagger}$ exists. Now,
 \begin{align*}
 D((T^{\dagger})^{\dagger}) &= R(T^{\dagger}) \oplus^{\perp} R(T^{\dagger})^{\perp}\\ 
 &= C(T) \oplus^{\perp} R(\overline{T^{\dagger}})^{\perp}\\
 &= C(T) \oplus^{\perp} N(({T^{\dagger}})^{*})\\
 &= C(T) \oplus^{\perp} N({(T^{*})^{\dagger}}) ~ (\text{Theorem \ref{thm 2.3}})\\
 & = C(T) \oplus^{\perp} N(\overline{T}) .
 \end{align*}
 So, The relations $N(\overline T) = \overline {N(T)}$ and $C(T) \subset C(\overline{T})$ confirm that $D(T) \subset D((T^{\dagger})^{\dagger}) \subset D(\overline{T})$. We know $\overline{T} = ((\overline{T})^{\dagger})^{\dagger}$. We claim that $((\overline{T})^{\dagger})^{\dagger} \supset (T^{\dagger})^{\dagger}$.
 \begin{align*}
  (T^{\dagger})^{\dagger}(x) &= \overline{T}(x) = 0, \text{ for all } x\in N(\overline{T})\\
  &= (T^{\dagger}\vert_{C(T^{\dagger})})^{-1}(x) = ({\overline{T}}^{\dagger}\vert_{C({\overline{T}}^{\dagger})})^{-1}(x)= T(x), \text{ for all } x \in C(T).
 \end{align*}
 Therefore, $T \subset (T^{\dagger})^{\dagger} \subset \overline{T}$.
 \end{proof}
 \begin{proposition}\label{pro 2.4}
 Let $T$ be a densely defined closable operator from $D(T) \subset H$ into $K$. Then $T^{\dagger}T$ and $TT^{\dagger}$ both are symmetric operators.
 \end{proposition}
 \begin{proof}
 We know that $T^{\dagger}T = P_{N(T)^{\perp}}\vert_{D(T)}$ and $TT^{\dagger} = P_{\overline{R(T)}}\vert_{D(T^{\dagger})}$, where $P_{N(T)^{\perp}}$ and $P_{\overline{R(T)}}$ both are orthogonal projection from $H$ onto $N(T)^{\perp}$ and $K$ onto $\overline{R(T)}$ respectively. Then, $(T^{\dagger}T)^{*} = P_{N(T)^{\perp}} \supset T^{\dagger}T$ and $(TT^{\dagger})^{*} = P_{\overline{R(T)}} \supset TT^{\dagger}$. Therefore, $T^{\dagger}T$ and $TT^{\dagger}$ both are symmetric operators. 
 \end{proof}
 \begin{theorem}\label{thm 2.5}
 Let $T$ be a closable operator from $D(T) \subset H$ into $K$. Then the following statements hold good:
 \begin{enumerate}
 \item $R(T^{\dagger}) = R(T^{\dagger}T)$ and $N(T) = N(T^{\dagger}T)$;
 \item $R(T) = R(TT^{\dagger})$ and $N(T^{\dagger}) = N(TT^{\dagger})$;
 \item When $T$ is densely defined then $\overline{R(T^{\dagger})} = \overline{C(T)} = N(T)^{\perp} = \overline{R(T^{*})}$.
 \end{enumerate}
 \end{theorem}
 \begin{proof}

     $(\mathit{1})$ $R(T^{\dagger}) = C(T) = R(T^{\dagger}T)$. 
     
     The inclusion $N(T) \subset N(T^{\dagger}T)$ is obvious. To prove the reverse inclusion, we consider an element $x \in N(T^{\dagger}T)$. So, $Tx \in R(T) \cap R(T)^{\perp}$ which implies $x \in N(T)$. Therefore, $N(T) = N(T^{\dagger}T)$.
     
     \noindent $(\mathit{2})$ The inclusion relation $R(TT^{\dagger}) \subset R(T)$ is easy to prove. To prove that reverse inclusion, let us consider an element  $Tx \in R(T)$. Then $Tx = TT^{\dagger}Tx \in R(TT^{\dagger})$. Therefore, $R(T) = R(TT^{\dagger})$.

     The inclusion relation $N(T^{\dagger}) \subset N(TT^{\dagger})$ holds. Now, we consider $y \in N(TT^{\dagger})$. Then $T^{\dagger}y \in C(T) \cap N(T) = \{0\}$ which confirms that $y\in N(T^{\dagger})$. Therefore, $N(T^{\dagger}) = N(TT^{\dagger})$.

     $(\mathit{3})$ By the definition of the Moore-Penrose inverse, we say $\overline{R(T^{\dagger})} = \overline{C(T)}$. Moreover, $\overline{C(T)} \subset N(T)^{\perp}$ is obvious. Let us consider an element $z \in N(T)^{\perp}$. There exists a sequence $\{z_{n}\}$ in $D(T)$ such that $z_{n} = z_{n}^{'} + z_{n}^{''} \to z$ as $n \to \infty$, where $z_{n}^{'} \in N(T)$ and $z_{n}^{''} \in C(T) \text{ for all } n \in \mathbb N$. So, $z_{n}^{''} \to z$ as $n \to \infty$ which implies that $z \in \overline{C(T)}$. Thus, $N(T)^{\perp} = \overline{C(T)}$.

   The relation $N(\overline{T}) = \overline{N(T)}$ guarantees that $N(T)^{\perp} = N(\overline{T})^{\perp} = \overline{R(\overline{T}^{*})} = \overline{R(T^{*})}$. 
   Therefore, the relation $\overline{R(T^{\dagger})} = \overline{C(T)} = N(T)^{\perp} = \overline{R(T^{*})}$ holds true.
  \end{proof}
  \begin{theorem}\label{thm 2.6}
  Let $T$ be a densely defined closable operator from $D(T) \subset H$ into $K$. Then $(T^{*}T)^{\dagger} \supset T^{\dagger}(T^{*})^{\dagger}$ (without assuming the condition $N(\overline{T}) = \overline{N(T)}$).
  \end{theorem}
  \begin{proof}
  Firstly, we claim that $T^{*}T$ has a decomposable domain. We know $D(T) = N(T) \oplus^{\perp} C(T)$ and $N(T^{*}T) = N(T)$. It is obvious to show that $N(T^{*}T) \oplus^{\perp} C(T^{*}T) \subset D(T^{*}T)$. Let us consider $x \in D(T^{*}T) \subset D(T)$. Then $x = x_{1} + x_{2}$, where $x_{1} \in N(T) = N(T^{*}T) \text { and } x_{2} \in C(T)$. So, $x_{2} \in D(T^{*}T)$. Again $x_{2} \in C(T) \cap D(T^{*}T) = N(T^{*}T)^{\perp} \cap D(T^{*}T) = C(T^{*}T)$ which implies that $x \in N(T^{*}T) \oplus^{\perp} C(T^{*}T)$. Thus, $N(T^{*}T) \oplus^{\perp} C(T^{*}T) = D(T^{*}T)$. The decomposable domain of $T^{*}T$ confirms that the existence of $(T^{*}T)^{\dagger}$.

  Now we will show that $D(T^{\dagger}(T^{*})^{\dagger}) \subset D((T^{*}T)^{\dagger})$. Let us consider an element $z \in D(T^{\dagger}(T^{*})^{\dagger})$. Then $z= z_{1} + z_{2} \in R(T^{*}) \oplus^{\perp} R(T^{*})^{\perp} \subset R(T^{*}) + R(T^{*}T)^{\perp}$. Moreover, $(T^{*})^{\dagger}z \in D(T^{\dagger}) = R(T) \oplus^{\perp} R(T)^{\perp}$. By the definition of the Moore-Penrose inverse, we get $(T^{*})^{\dagger}z \in C(T^{*}) = N(T^{*})^{\perp} \cap D(T^{*}) = \overline{R(T)} \cap D(T^{*})$. So, $(T^{*})^{\dagger}z \in R(T) \cap D(T^{*})$. We get $w_{1} \in H$ such that 
  \begin{align*}
  (T^{*})^{\dagger}z = Tw_{1}\\
  (T^{*})^{\dagger}z_{1} = Tw_{1}\\
  z_{1} = T^{*}Tw_{1}.
\end{align*}
This confirms that $z = z_{1} + z_{2} \in R(T^{*}T) \oplus^{\perp} R(T^{*}T)^{\perp} = D(T^{*}T)^{\dagger}$. Hence, $D(T^{\dagger}(T^{*})^{\dagger}) \subset D((T^{*}T)^{\dagger})$. To prove the required relation, we consider $s\in D(T^{\dagger} (T^{*})^{\dagger}) \cap R(T^{*})^{\perp}$. Since, $R(T^{*})^{\perp} \subset R(T^{*}T)^{\perp}$ which implies $T^{\dagger}(T^{*})^{\dagger}s = (T^{*}T)^{\dagger}s = 0$. Let us look an element $u \in D(T^{\dagger}(T^{*})^{\dagger}) \cap R(T^{*})$. There exists an element $v \in C(T^{*}) = \overline{R(T)} \cap D(T^{*})$ such that $T^{*}v = u$. Furthermore, $(T^{*})^{\dagger}u = v \in D(T^{\dagger}) = R(T) \oplus^{\perp} R(T)^{\perp}$. So, $v \in R(T) \cap D(T^{*})$ which guarantees that the existence of an element $z_{0} \in C(T)$ such that $Tz_{0} = v$. Again $T^{\dagger}v= z_{0}$ implies $T^{\dagger}(T^{*})^{\dagger}u = z_{0}$. It is clear that $z_{0} \in D(T^{*}T) \cap N(T)^{\perp} = D(T^{*}T) \cap N(T^{*}T)^{\perp} = C(T^{*}T)$. Then, $(T^{*}T)^{\dagger}u = (T^{*}T)^{\dagger} (T^{*}T)z_{0} = z_{0}$. Therefore, $T^{\dagger} (T^{*})^{\dagger} \subset (T^{*}T)^{\dagger}$.
\end{proof}
\begin{theorem}\label{thm 2.7}
Let $T$ be a densely defined closable operator from $D(T) \subset H$ into $K$. Then $(TT^{*})^{\dagger} \supset (T^{*})^{\dagger} T^{\dagger}$.
\end{theorem}
\begin{proof}
Firstly, we claim that $D(TT^{*})$ is decomposable. It is obvious that $N(TT^{*}) \oplus^{\perp} C(TT^{*}) \subset D(TT^{*})$ and $N(TT^{*}) = N(T^{*})$. To prove the reverse inclusion, we consider $x \in D(TT^{*}) \subset D(T^{*}) = N(TT^{*}) \oplus^{\perp} C(T^{*})$. Now, we can write $x = x_{1} + x_{2}$, where $x_{1} \in N(TT^{*})$ and $x_{2} \in C(T^{*})$. So, $x_{2} \in C(T^{*}) \cap D(TT^{*}) = C(TT^{*})$ which shows that $x \in N(TT^{*}) \oplus^{\perp} C(TT^{*})$. Thus, $D(TT^{*})$ is a decomposable domain of $TT^{*}$. Moreover, $(TT^{*})^{\dagger}$ exists. 

Now, we will show that $D((T^{*})^{\dagger}T^{\dagger}) \subset D((TT^{*})^{\dagger})$. Let us consider an element $y \in D((T^{*})^{\dagger}T^{\dagger})$. Then $y = y_{1} + y_{2} \in R(T) \oplus^{\perp} R(T)^{\perp} \subset R(T) + R(TT^{*})^{\perp}$, where $y_{1} \in R(T)$ and $y_{2} \in R(T)^{\perp}$. Again, $T^{\dagger}y_{1} \in R(T^{*}) \oplus^{\perp} R(T^{*})^{\perp}$ and $T^{\dagger}y_{1} \in C(T) = N(\overline{T})^{\perp} \cap D(T) = \overline{R(T^{*})} \cap D(T)$ which implies $T^{\dagger}y_{1} \in R(T^{*}) \cap D(T)$. There exists $w \in C(T^{*})$ such that 
\begin{align*}
T^{\dagger}y_{1} = T^{*}w\\
y_{1}= TT^{*}w.
\end{align*}
So, $y \in R(TT^{*}) \oplus^{\perp} R(TT^{*})^{\perp} = D((TT^{*})^{\dagger})$ which implies that $D((T^{*})^{\dagger}T^{\dagger}) \subset D(TT^{*})^{\dagger}$.
To establish the required statement, we consider $s \in D((T^{*})^{\dagger}T^{\dagger}) \cap R(T)^{\perp}$. Since $R(T)^{\perp} \subset R(TT^{*})^{\perp}$. Thus, $(T^{*})^{\dagger} T^{\dagger}s = (TT^{*})^{\dagger}s= 0$. Now, let us look an element $z \in D((T^{*})^{\dagger}T^{\dagger}) \cap R(T)$. There exists $z_{0} \in C(T)= N(\overline{T})^{\perp} \cap D(T)= \overline{R(T^{*})} \cap D(T)$ such that $Tz_{0} = z$. Again, $T^{\dagger}z= z_{0} \in D((T^{*})^{\dagger}) = R(T^{*}) \oplus^{\perp} R(T^{*})^{\perp}$ which implies that $z_{0} \in R(T^{*}) \cap D(T)$. Now, we have an element $v \in C(T^{*}) = N(T^{*})^{\perp} \cap D(T^{*})$ such that $T^{*}v= z_{0}$. Furthermore, $(T^{*})^{\dagger}T^{\dagger}z= v$. It is easy to show that $v \in D(TT^{*}) \cap N(TT^{*})^{\perp} = C(TT^{*})$. Hence, $(TT^{*})^{\dagger}z = (TT^{*})^{\dagger} (TT^{*})v= v$. Therefore, $(T^{*})^{\dagger}T^{\dagger} \subset (TT^{*})^{\dagger}$.
\end{proof}
\begin{corollary}\label{cor 2.8}
Let $T$ be a densely defined closable operator from $D(T) \subset H$ into $K$ with $R(T^{*})$ be closed. Then $T^{*}(TT^{*})^{\dagger} \supset T^{\dagger}$.
\end{corollary}
\begin{proof}
From Theorem \ref{thm 2.7}, we say that $T^{*}(TT^{*})^{\dagger} \supset T^{*}(T^{*})^{\dagger}T^{\dagger} = (P_{R(T^{*})}\vert_{D((T^{*})^{\dagger})})T^{\dagger}$. Now, $R(T^{\dagger}) = C(T)= N(\overline{T})^{\perp} \cap D(T)= R(T^{*}) \cap D(T)$. Therefore, $T^{*}(TT^{*})^{\dagger} \supset T^{\dagger}$.
\end{proof}
\begin{corollary}\label{cor 2.9}
Let $T$ be a densely defined closable operator from $D(T) \subset H$ into $K$ with $R(T)$ be closed. Then $(T^{*}T)^{\dagger}T^{*} \subset T^{\dagger}$.
\end{corollary}
\begin{proof}
 From Theorem \ref{thm 2.6}, we get $(T^{*}T)^{\dagger}T^{*} \supset T^{\dagger}(T^{*})^{\dagger}T^{*}$. Since $R(T)$ is closed. So, $R(T^{*})^{\dagger} = C(T^{*}) = N(T^{*})^{\perp} \cap D(T^{*})= R(T) \cap D(T^{*})$. Thus, $D(T^{\dagger}(T^{*})^{\dagger} T^{*}) = D(T^{*}) \supset D((T^{*}T)^{\dagger}T^{*})$ which implies $(T^{*}T)^{\dagger}T^{*} = T^{\dagger}(T^{*})^{\dagger}T^{*}$. Let us consider an element $x \in N(T^{*})= R(T)^{\perp}$. Then $T^{\dagger}x = T^{\dagger}(T^{*})^{\dagger}T^{*}x =0$. Similarly, for an element $z \in N(T^{*})^{\perp} \cap D(T^{*}) \subset R(T)$, we have $T^{\dagger}(T^{*})^{\dagger}T^{*} z= T^{\dagger}z$. Therefore, $(T^{*}T)^{\dagger}T^{*} \subset T^{\dagger}$.
 \end{proof}
 \begin{theorem}\label{thm 2.10}
 Let $T$ be a closable operator from $D(T) \subset H$ into $K$. Then the following statements hold good:
 \begin{enumerate}
     \item $D(TT^{\dagger}) = D(T^{\dagger})$ and $D(T^{\dagger}T) = D(T)$;
     \item $R(TT^{\dagger}) = R(T) \text{ and } R(T^{\dagger}T) = R(T^{\dagger})$;
     \item when $T$ is densely defined then $N(\overline{T}) = N((T^{*})^{\dagger})$.
 \end{enumerate}
 \end{theorem}
 \begin{proof}
 $(\mathit{1})$ $D(TT^{\dagger}) \subset D(T^{\dagger})$ and $R(T^{\dagger}) = C(T) \subset D(T)$ confirm that $D(TT^{\dagger}) = D(T^{\dagger})$.

 Moreover, $D(T^{\dagger}T) \subset D(T)$ and $R(T) \subset D(T^{\dagger})$ guarantee that $D(T^{\dagger}T) = D(T)$.
 
 $(\mathit{2})$ It is obvious that $R(TT^{\dagger}) \subset R(T)$. To prove the reverse inclusion, let us consider an element $y\in R(T)$. Then there exist $x \in C(T)$ such that $Tx =y$ and $T^{\dagger}y = x$ which implies $TT^{\dagger}y = y \in R(TT^{\dagger})$. Hence, $R(T) = R(TT^{\dagger})$.

 Furthermore, $R(T^{\dagger}T) \subset R(T^{\dagger})$. Now, consider an element $x \in R(T^{\dagger})$. There exists an element $z \in R(T)$ such that  $T^{\dagger}z = x$. Again, we have an element $w \in C(T)$ with $Tw =z$. Thus, $x = T^{\dagger}Tw \in R(T^{\dagger}T)$. Therefore, $R(T^{\dagger}T) = R(T^{\dagger})$.

  $(\mathit{3})$ $N((T^{*})^{\dagger}) = R(T^{*})^{\perp} = N(T^{**})= N(\overline{T})$.
 \end{proof}
 \begin{theorem}\label{thm 2.11}
 Let $T$ be a densely defined closable operator from $D(T) \subset H$ into $K$ with $R(T^{*})$ be closed and $R(T^{*}) \subset D(T)$. Then $T \subset TT^{*}(T^{*})^{\dagger}$. 
 \end{theorem}
 \begin{proof}
 Firstly, we claim that $D(T) \subset D(TT^{*}(T^{*})^{\dagger})$. Let $x \in D(T)$. Then $x \in C(T) \oplus^{\perp} N(T) = (N(\overline{T})^{\perp} \cap D(T)) \oplus^{\perp} N(T) = R(T^{*}) \oplus^{\perp} N(T) \subset D((T^{*})^{\dagger})$. The mentioned two relations $R(T^{*}) \subset D(T)$ and $R((T^{*})^{\dagger})= C(T^{*}) \subset D(T^{*})$ confirm that $x \in D(TT^{*}(T^{*})^{\dagger})$. So, $D(T) \subset D(TT^{*}(T^{*})^{\dagger})$. The closeness of $R(T^{*})$ says that $R(T^{*}) = N(T)^{\perp}$. Let us consider an element $z = z_{1} + z_{2} \in N(T) \oplus^{\perp} C(T)$, where $z_{1} \in N(T) \text{ and } z_{2}\in C(T)$. Then, $TT^{*}(T^{*})^{\dagger} z =  Tz_{2} = Tz$. Therefore, $T \subset TT^{*}(T^{*})^{\dagger}$.
 \end{proof}
 \begin{theorem}\label{thm 2.12}
 Let $T$ be a densely defined closable operator from $D(T) \subset H$ into $K$ with the condition $R(T) \subset D(T^{*})$. Then $T = (T^{*})^{\dagger}T^{*}T$.
 \end{theorem}
 \begin{proof}
 It is obvious that $D(T) \supset D((T^{*})^{\dagger}T^{*}T)$. From the given condition $R(T) \subset D(T^{*})$ says that the inclusion $ D(T) \subset D((T^{*})^{\dagger}T^{*}T)$ which implies $D(T) = D((T^{*})^{\dagger}T^{*}T)$. Now let us consider an element $x \in D(T)$. Then, $(T^{*})^{\dagger}T^{*}Tx= P_{\overline{R(T)}}\vert_{D(T^{*})}(Tx) = Tx$. Therefore, $T = (T^{*})^{\dagger}T^{*}T$.
 \end{proof}
 \begin{theorem}\label{thm 2.13}
 Let $T$ be a densely defined closable operator from $D(T) \subset H$ into $K$ with $R(T)$ be closed. Then $(T^{\dagger})^{*} = TT^{\dagger}(T^{\dagger})^{*}$.
 \end{theorem}
 \begin{proof}
 It is obvious to prove the inclusion $D(TT^{\dagger}(T^{\dagger})^{*}) \subset D((T^{\dagger})^{*})$. From Theorem \ref{thm 2.2}, we get $R((T^{\dagger})^{*}) = R((T^{*})^{\dagger}) = C(T^{*}) = N(T^{*})^{\perp} \cap D(T^{*}) = R(T) \cap D(T^{*})$ which confirms that $D(TT^{\dagger}(T^{\dagger})^{*}) \supset D((T^{\dagger})^{*})$. Thus, $D(TT^{\dagger}(T^{\dagger})^{*}) = D((T^{\dagger})^{*})$. Now, consider an element $x \in D(TT^{\dagger}(T^{\dagger})^{*})$, we get $x = x_{1} + x_{2} \in R(T^{*}) \oplus^{\perp} R(T^{*})^{\perp}$, where $x_{1} \in R(T^{*})$ and $x_{2} \in R(T^{*})^{\perp}$. Then there exists an element $w \in C(T^{*}) \subset R(T)$ such that $x_{1} = T^{*}w$. Moreover,
 \begin{equation}
 TT^{\dagger}(T^{\dagger})^{*}x= TT^{\dagger}w= w =(T^{*})^{\dagger}T^{*}w= (T^{*})^{\dagger}x_{1}= (T^{*})^{\dagger}x.
 \end{equation}
 Therefore, $(T^{\dagger})^{*} = TT^{\dagger}(T^{\dagger})^{*}$.
 \end{proof}
 \begin{theorem}\label{thm 2.14}
 Let $T$ be a densely defined closable operator from $D(T) \subset H$ into $K$ with $R(T^{*})$ be closed and $R(T^{*}) \subset D(T)$. Then $(T^{\dagger})^{*}T^{\dagger}T \subset (T^{\dagger})^{*}$.
 \end{theorem}
 \begin{proof}
 Since, $R(T^{*})$ is closed. Then $D((T^{\dagger})^{*}) = D((T^{*})^{\dagger}) = H \supset D((T^{\dagger})^{*}T^{\dagger}T)$. Let us consider an element $x \in D((T^{\dagger})^{*}T^{\dagger}T)$. So, $x= x_{1} + x_{2} \in H$, where $x_{1}\in R(T^{*})$ and $x_{2} \in R(T^{*})^{\perp}$. From Theorem \ref{thm 2.2}, we get $(T^{\dagger})^{*}T^{\dagger}Tx = (T^{\dagger})^{*}P_{R(T^{*})}\vert_{D(T)} (x)= (T^{\dagger})^{*}x_{1} = (T^{\dagger})^{*}x$. Therefore, $(T^{\dagger})^{*}T^{\dagger}T \subset (T^{\dagger})^{*}$.
 \end{proof}
 \begin{theorem}\label{thm 2.15}
 Let $T$ be a densely defined closable operator from $D(T) \subset H$ into $K$. Then $T^{*} = (T^{\dagger}T)^{*}T^{*}$.
 \end{theorem}
 \begin{proof}
 The denseness of $D(T)$ and $D(T^{\dagger}T) = D(T)$ confirm that the existence of $(T^{\dagger}T)^{*}$. It is obvious that $D(T^{*}) \supset D((T^{\dagger}T)^{*}T^{*})$. Now, we will show the reverse inclusion. It is enough to prove that $R(T^{*}) \subset D((T^{\dagger}T)^{*})$. The relations $(T^{\dagger}T)^{*} \supset T^{*}(T^{\dagger})^{*} = T^{*}(T^{*})^{\dagger} $ and $R((T^{*})^{\dagger}) \subset D(T^{*})$ say that 
 \begin{align*}
 D((T^{\dagger}T)^{*}) \supset D(T^{*} (T^{*})^{\dagger}) = D((T^{*})^{\dagger}) \supset R(T^{*}).
 \end{align*}
 Thus, $D(T^{*}) = D((T^{\dagger}T)^{*}T^{*})$. Now consider an element $y \in D((T^{\dagger}T)^{*}T^{*})$. Then, 
 \begin{align*}
 (T^{\dagger}T)^{*}T^{*}y &= (P_{N(T)^{\perp}}\vert_{D(T)})^{*}T^{*}y\\
 &= P_{N(T)^{\perp}} T^{*}y\\
 &= P_{N(\overline{T})^{\perp}} T^{*}y \\
 &= P_{\overline{R(T^{*})}}T^{*}y \\
 &= T^{*}y.
 \end{align*}
 Therefore, $T^{*} = (T^{\dagger}T)^{*}T^{*}$.
 \end{proof}
 \begin{theorem}\label{thm 2.16}
 Let $T$ be a densely defined closable operator from $D(T) \subset H$ into $K$. Then $T^{*} \supset T^{\dagger}TT^{*}$.
 \end{theorem}
 \begin{proof}
 We know, $D(T^{\dagger}TT^{*}) \subset D(T^{*})$. Let us consider an element $y \in D(T^{\dagger}TT^{*})$. Then $T^{*}y \in R(T^{*}) \cap D(T) \subset N(\overline{T})^{\perp} \cap D(T) = N(T)^{\perp} \cap D(T)$. Moreover, $T^{\dagger}TT^{*}y = P_{N(T)^{\perp}}\vert_{D(T)} T^{*}y= T^{*}y$. Therefore, $T^{*} \supset T^{\dagger}TT^{*}$.
 \end{proof}
 Now, we will present an illustrative example to justify our results.
 \begin{example}
	 Let $T_{\phi}$ be a multiplication operator on $L^{2}(M)$, where $L^{2}(M) = L^{2} (M, \Sigma, \mu)$, $\mu$ is a Borel regular measure, $M \subset \mathbb R$  and $\mu(M) < \infty$. Define $T_{\phi}$ as $T_{\phi}(f) =\phi f$, for $f\in D(T_{\phi}) = \{f\in L^{2}(M) : \phi f \in L^{2}(M), \phi \in \mathcal{C} (M) , |\phi(x)| \geq 1\}$. Here, $\mathcal{C}(M)$ is the set of all continuous functions from $M$ to $\mathbb{C}$. Then $T_{\phi}$ is closed and $T_{\phi} = T_{\phi}^{*} =T_{\overline {\phi}}$, where $\overline{\phi}$ is the complex conjugate  of $\phi$.  Moreover, $R(T_{\phi}) =R{(T_{\phi}}^{*}) = L^{2}(M)$ \cite{majumdar2024hyers} and $N(T_{\phi}) = \overline {N(T_{\phi}|_{{\mathcal{C}}_{0}^{\infty}})} = \{0\}$, where ${T_{\phi}} |_{{\mathcal{C}}_{0}^{\infty}}$ is the restriction of $T_{\phi}$ in the domain ${\mathcal{C}}_{0}^{\infty}$. Again, ${T_{\phi}} |_{{\mathcal{C}}_{0}^{\infty}}$ is not closed but it is closable because  ${T_{\phi}} |_{{\mathcal{C}}_{0}^{\infty}} \subsetneq \overline{{T_{\phi}} |_{{\mathcal{C}}_{0}^{\infty}}} = T_{\phi}$. ${T_{\phi}} |_{{\mathcal{C}}_{0}^{\infty}}$ has decomposable domain because $D({T_{\phi}} |_{{\mathcal{C}}_{0}^{\infty}})= C({T_{\phi}} |_{{\mathcal{C}}_{0}^{\infty}})$. Thus, $({T_{\phi}} |_{{\mathcal{C}}_{0}^{\infty}})^{\dagger}$ exists and $D(({T_{\phi}} |_{{\mathcal{C}}_{0}^{\infty}})^{\dagger}) = R({T_{\phi}} |_{{\mathcal{C}}_{0}^{\infty}})$. We also get $(({T_{\phi}} |_{{\mathcal{C}}_{0}^{\infty}})^{\dagger})^{\dagger} = ({T_{\phi}} |_{{\mathcal{C}}_{0}^{\infty}}) \subsetneq T_{\phi}$. 
  It is obvious to show that $R(T_{\overline{\phi}}) = L^{2}(M)= R(T_{\phi}^{*}) = R(({T_{\phi}} |_{{\mathcal{C}}_{0}^{\infty}})^{*})$. $D(({T_{\phi}} |_{{\mathcal{C}}_{0}^{\infty}})^{\dagger})$ is dense in $L^{2}(M)$ which implies that the existence of $(({T_{\phi}} |_{{\mathcal{C}}_{0}^{\infty}})^{\dagger})^{*}$. Theorem \ref{thm 2.2} says that $(({T_{\phi}} |_{{\mathcal{C}}_{0}^{\infty}})^{\dagger})^{*} = (T_{\overline{\phi}})^{\dagger}$ with the whole domain $L^{2}(M)$. Therefore, $g_{1} \in L^{2}(M)$, we have
  \begin{align*}
  (({T_{\phi}} |_{{\mathcal{C}}_{0}^{\infty}})^{\dagger})^{*}(g_{1}) = g_{2}, \text{ where } g_{2}\in D(T_{\overline{\phi}}) \text{ and } T_{\overline{\phi}}(g_{2})= g_{1}.
  \end{align*}
\end{example}
\section{Properties of Moore-Penrose inverses of the direct sum of closed operators in Hilbert spaces:}
Let $H$ and $K$ be two Hilbert spaces. The space $H \bigoplus K$ defined by $H \bigoplus K = \{(h, k): h\in H, k\in K\}$ is a linear space with respect to addition and scalar multiplication defined by 
\begin{align*}
&~~~~~~~~~~~~~~~~~~~~~~(h_{1}, k_{1}) + (h_{2}, k_{2}) = (h_{1} + h_{2}, k_{1} + k_{2}), \text{ and }\\
&\lambda (h, k ) = (\lambda h, \lambda k), \text{ for all } h, h_{1}, h_{2}\in H, \text{ for all } k, k_{1}, k_{2}\in K \text{ and } \lambda\in \mathbb{K}, ~(\mathbb{K} = \mathbb{R} \text{ or } \mathbb{C}).
\end{align*}
Now, $H \bigoplus K$ is an inner product space with respect to the inner product given by
\begin{align*}
\langle (h_{1}, k_{1}) (h_{2}, k_{2}) \rangle = \langle h_{1}, h_{2}\rangle + \langle k_{1}, k_{2}\rangle, \text{ for all } h_{1}, h_{2}\in H, \text{ and } \text{ for all }k_{1}, k_{2} \in K. 
\end{align*}
The norm on $H \bigoplus K$ is defined by 
\begin{align*}
\|(h, k)\| = (\|h\|^{2} + \|k\|^{2})^{\frac{1}{2}}, \text{ for all } (h,k) \in H \times K.
\end{align*}
Moreover, The direct sum of two operators $T_{1} \text{ and } T_{2}$ from $D(T_{1}) \subset H_{1}$ to $K_{1}$ and from $D(T_{2}) \subset H_{2}$ to $K_{2}$ respectively is defined by 
\begin{align*}
(T_{1} \bigoplus T_{2}) (h_{1}, h_{2}) = (T_{1}h_{1}, T_{2}h_{2}), \text{ for all } h_{1} \in D(T_{1}), \text { and for all } h_{2}\in D(T_{2}).
\end{align*}
\begin{theorem}\label{thm 3.1}
Let $T_{1}: D(T_{1}) \subset H_{1} \to K_{1}$ and $T_{2}: D(T_{2}) \subset  H_{2} \to K_{2}$ be two closed operators with closed ranges. Then $T = T_{1} \bigoplus T_{2} : D(T_{1}) \bigoplus D(T_{2}) \subset H_{1} \bigoplus H_{2} \to K_{1} \bigoplus K_{2}$ has the Moore-Penrose inverse. Moreover,
 \begin{align*}
 T^{\dagger} = (T_{1} \bigoplus T_{2})^{\dagger} = T_{1}^{\dagger} \bigoplus T_{2}^{\dagger}.
 \end{align*}
\end{theorem}
\begin{proof}
$T$ is closed because $T_{1}$ and $T_{2}$ both are closed. So, the Moore-Penrose inverse $T^{\dagger}$ exists in domain $D(T^{\dagger}) = R(T) \oplus R(T)^{\perp}$. It is obvious to show that $R(T) = R(T_{1}) \bigoplus R(T_{2})$ and $N(T) = N(T_{1}) \bigoplus N(T_{2})$. Here, the closed property of $R(T_{1})$ and $R(T_{2})$ guarantees that $R(T)$ is also closed. We will show that 
\begin{align*}
 P_{R(T_{1}) \bigoplus R(T_{2})} = P_{R(T_{1})} \bigoplus P_{R(T_{2})}.
\end{align*}
For all $h_{1} \in D(T_{1}) \text{ and } h_{2} \in D(T_{2})$, 
\begin{align*}
 P_{R(T_{1}) \bigoplus R(T_{2})}(T_{1}h_{1}, T_{2}h_{2}) &= (T_{1}h_{1}, T_{2}h_{2})\\
 &= (P_{R(T_{1})}(T_{1}h_{1}) , P_{R(T_{2})}(T_{2}h_{2}))\\
 &= P_{R(T_{1})} \bigoplus P_{R(T_{2})} (T_{1}h_{1}, T_{2}h_{2}).
 \end{align*}
 Now, let us consider $(k_{1}, k_{2}) \in (R(T_{1}) \bigoplus R(T_{2}))^{\perp}$. Since $R(T_{1}) \text{ and } R(T_{2})$ both are closed then there exist $k_{1}^{'}, k_{1}^{''}, k_{2}^{'}, k_{2}^{''}$
in $R(T_{1}), R(T_{1})^{\perp}, R(T_{2}) \text{ and } R(T_{2})^{\perp}$ respectively such that $k_{1} = k_{1}^{'} + k_{1}^{''}$ and $k_{2} = k_{2}^{'} + k_{2}^{''}$. Again, for all $(h_{1}, h_{2}) \in D(T_{1}) \bigoplus D(T_{2})$, 
\begin{equation}\label{eqn 4}
 \langle T_{1}h_{1}, k_{1}^{'}\rangle + \langle T_{2}h_{2}, k_{2}^{'}\rangle = \langle  (T_{1}h_{1}, T_{2}h_{2}), (k_{1}, k_{2}) \rangle = 0.
\end{equation}
there exist $v_{1} \in D(T_{1})$ and $v_{2} \in D(T_{2})$ such that $T_{1}v_{1} = k_{1}^{'}$ and $T_{2}v_{2} = k_{2}^{'}$. We consider $h_{1} = v_{1}$ and $h_{2}= 0$ in (\ref{eqn 4}), we get $k_{1}^{'} =0$. Similarly, when $h_{1}=0 \text{ and } h_{2} = v_{2}$ then $k_{2}^{'} = 0$. Thus, $(k_{1}, k_{2}) = (k_{1}^{''}, k_{2}^{''}) \in R(T_{1})^{\perp} \bigoplus R(T_{2})^{\perp}$. So,
\begin{align*}
P_{R(T_{1}) \bigoplus R(T_{2})}(k_{1}, k_{2})= (0,0)
\end{align*}
and
\begin{align*}
P_{R(T_{1})} \bigoplus P_{R(T_{2})}(k_{1}, k_{2}) & = P_{R(T_{1})} \bigoplus P_{R(T_{2})}(k_{1}^{''}, k_{2}^{''})\\
&= (P_{R(T_{1})}k_{1}^{''} , P_{R(T_{2})}k_{2}^{''})\\
&= (0,0).
\end{align*}
Our above claim is proven which implies that
\begin{equation}\label{eqn 5}
TT^{\dagger} = P_{R(T)} = P_{R(T_{1}) \bigoplus R(T_{2})} 
= P_{R(T_{1})} \bigoplus P_{R(T_{2})}
= T_{1}T_{1}^{\dagger} \bigoplus T_{2} T_{2}^{\dagger}
= T (T_{1}^{\dagger} \bigoplus T_{2}^{\dagger}).
\end{equation}
We claim that $P_{(N(T_{1}) \bigoplus N(T_{2}))^{\perp}}\vert_{ D(T)} = P_{N(T_{1})^{\perp}} \vert_ {D(T_{1})} \bigoplus P_{N(T_{2})^{\perp}}\vert_{ D(T_{2})}$.
Let us consider an arbitrary element $(x_{1}, x_{2}) \in N(T) = N(T_{1}) \bigoplus N(T_{2}) \subset D(T)$. Then,
$P_{(N(T_{1}) \bigoplus N(T_{2}))^{\perp}}\vert_{ D(T)}(x_{1}, x_{2}) = (P_{N(T_{1})^{\perp}} \vert_ {D(T_{1})}(x_{1}) , P_{N(T_{2})^{\perp}}\vert_{ D(T_{2})} (x_{2})) = (0,0)$. Again, for an element $(w_{1}, w_{2}) \in (N(T_{1}) \bigoplus N(T_{2}))^{\perp} \cap D(T)$, we have 
\begin{align*}
P_{(N(T_{1}) \bigoplus N(T_{2}))^{\perp}}\vert_{ D(T)}(w_{1}, w_{2}) = (w_{1}, w_{2}).
\end{align*}
For all $(s_{1}, s_{2}) \in N(T_{1}) \bigoplus N(T_{2})$, we get 
\begin{equation}\label{eqn 6}
\langle s_{1}, w_{1}\rangle + \langle s_{2}, w_{2} \rangle = \langle (s_{1}, s_{2}) , (w_{1}, w_{2})\rangle = 0.
 \end{equation}
 There are elements $w_{1}^{'}, w_{1}^{''}, w_{2}^{'} \text{ and } w_{2}^{''}$ in $N(T_{1}), C(T_{1}), N(T_{2}) \text{ and } C(T_{2})$ respectively such that $w_{1} = w_{1}^{'} + w_{1}^{''} \text{ and } w_{2}= w_{2}^{'} + w_{2}^{''}$. From the equation (\ref{eqn 6}), we get
 \begin{equation}\label{eqn 7}
 \langle s_{1}, w_{1}^{'}\rangle + \langle s_{2}, w_{2}^{'} \rangle = \langle (s_{1}, s_{2}) , (w_{1}, w_{2}) \rangle = (0,0).
 \end{equation}
 Let us take $s_{1} = w_{1}^{'} \text { and } s_{2}=0$, we have $w_{1}^{'} =0$. Similarly, $s_{1} = 0 \text{ and } s_{2} = w_{2}^{'}$ say that $w_{2}^{'} = 0$. Thus, $(w_{1}, w_{2}) = (w_{1}^{''}, w_{2}^{''})$ which implies 
 \begin{align*}
( P_{N(T_{1})^{\perp}} \vert_ {D(T_{1})} \bigoplus P_{N(T_{2})^{\perp}}\vert_{ D(T_{2})} )(w_{1}, w_{2}) &= (P_{N(T_{1})^{\perp}} \vert_ {D(T_{1})}(w_{1}^{''}), P_{N(T_{2})^{\perp}}\vert_{ D(T_{2})}(w_{2}^{''}))\\
&= (w_{1}, w_{2}).
 \end{align*}
 Hence, the relation $P_{(N(T_{1}) \bigoplus N(T_{2}))^{\perp}}\vert_{ D(T)} = P_{N(T_{1})^{\perp}} \vert_ {D(T_{1})} \bigoplus P_{N(T_{2})^{\perp}}\vert_{ D(T_{2})}$ is true. So,
 \begin{equation}\label{eqn 8}
 T^{\dagger}T = P_{(N(T_{1}) \bigoplus N(T_{2}))^{\perp}}\vert_{ D(T)}= P_{N(T_{1})^{\perp}} \vert_ {D(T_{1})} \bigoplus P_{N(T_{2})^{\perp}}\vert_{ D(T_{2})} = T_{1}^{\dagger}T_{1} \bigoplus T_{2}^{\dagger}T_{2}= (T_{1}^{\dagger} \bigoplus T_{2}^{\dagger}) T. 
 \end{equation}
 Furthermore,
 \begin{equation}\label{eqn 9}
 (T_{1}^{\dagger} \bigoplus T_{2}^{\dagger}) (T_{1} \bigoplus T_{2}) (T_{1}^{\dagger} \bigoplus T_{2}^{\dagger}) = (T_{1}^{\dagger}T_{1}T_{1}^{\dagger} \bigoplus T_{2}^{\dagger}T_{2}T_{2}^{\dagger})=(T_{1}^{\dagger} \bigoplus T_{2}^{\dagger}).
 \end{equation}
 Therefore, Theorem 5.7 \cite{MR0451661} and the equations (\ref{eqn 5}), (\ref{eqn 8}) and (\ref{eqn 9}) justify the relation $T^{\dagger} = T_{1}^{\dagger} \bigoplus T_{2}^{\dagger}$.
 \end{proof}
 \begin{corollary}
 Let $T_{i} : D(T_{i}) \subset H_{i} \to K_{i} ~(i = 1, 2 , \dots, n)$ be closed operators with closed ranges $R(T_{i})~ (i = 1, 2, \dots, n$). Then the Moore-Penrose inverse of $T_{1} \bigoplus T_{2} \bigoplus \dots \bigoplus T_{n}$ exists. Moreover, $(T_{1} \bigoplus T_{2} \bigoplus \dots \bigoplus T_{n})^{\dagger} = T_{1}^{\dagger} \bigoplus T_{2}^{\dagger} \bigoplus \dots \bigoplus T_{n}^{\dagger}$.
 \end{corollary}
 \begin{proof}
 By induction hypothesis and Theorem \ref{thm 3.1} prove the relation 
 \begin{align*}
 (T_{1} \bigoplus T_{2} \bigoplus \dots \bigoplus T_{n})^{\dagger} = T_{1}^{\dagger} \bigoplus T_{2}^{\dagger} \bigoplus \dots \bigoplus T_{n}^{\dagger}.
 \end{align*}
\end{proof}
\begin{lemma}\label{lemma 3.3}
Let $T_{1}: D(T_{1}) \subset H_{1} \to K_{1}$ and $T_{2}: D(T_{2}) \subset H_{2} \to K_{2}$ be two densely defined operators. Then $T_{1} \bigoplus T_{2}$ is a densely defined operator with
\begin{align*}
(T_{1} \bigoplus T_{2})^{*} = T_{1}^{*} \bigoplus T_{2}^{*}.
\end{align*}
\end{lemma}
\begin{proof}
 It is obvious that $D(T_{1} \bigoplus T_{2}) = D(T_{1}) \bigoplus D(T_{2})$ is densely defined in $H_{1} \bigoplus H_{2}$. So, $(T_{1} \bigoplus T_{2})^{*}$ exists. 
 \begin{align*}
 D(T_{1} \bigoplus T_{2})^{*} &= \{ (u,v) : \langle (T_{1} \bigoplus T_{2})(h_{1}, h_{2}), (u,v)\rangle \text{ is bounded for all} ~(h_{1}, h_{2}) \in D(T_{1} \bigoplus T_{2})\}\\
 &= \{ (u,v) : \langle T_{1}h_{1}, u\rangle + \langle T_{2}h_{2}, v \rangle \text{ is bounded for all} ~(h_{1}, h_{2}) \in D(T_{1} \bigoplus T_{2})\}\}
 \end{align*}
 Now consider, $h_{2} =0$. From above relation we get $u \in D(T_{1}^{*})$. Similarly, $h_{1} =0$ shows that $v \in D(T_{2}^{*})$. Thus, 
 \begin{equation}\label{eqn 10}
 D(T_{1} \bigoplus T_{2})^{*} \subset D(T_{1}^{*}) \bigoplus D(T_{2}^{*}).
 \end{equation}
 Taking an arbitrary element $(w_{1},w_{2}) \in D(T_{1}^{*}) \bigoplus D(T_{2}^{*})$, for all $(h_{1} , h_{2}) \in D(T_{1}) \bigoplus D(T_{2})$, we get
 \begin{align*}
 |\langle (T_{1} \bigoplus T_{2})(h_{1}, h_{2}) , (w_{1}, w_{2}) \rangle| &= |\langle T_{1}h_{1}, w_{1}\rangle + \langle T_{2}h_{2}, w_{2}\rangle|\\
 &= |\langle (h_{1}, h_{2}) (T_{1}^{*}w_{1}, T_{2}^{*} w_{2}) \rangle|\\
 & \leq ||(h_{1},h_{2})|| ||(T_{1}^{*}w_{1}, T_{2}^{*} w_{2})||
 \end{align*}
 Hence, $D(T_{1}^{*}) \bigoplus D(T_{2}^{*}) \subset  D(T_{1} \bigoplus T_{2})^{*}$ with 
 \begin{equation}\label{eqn 11}
 (T_{1} \bigoplus T_{2})^{*}(w_{1}, w_{2}) = (T_{1}^{*}w_{1}, T_{2}^{*}w_{2}).
 \end{equation}
 By the equations (\ref{eqn 10}) and (\ref{eqn 11}) guarantee that $(T_{1} \bigoplus T_{2})^{*} = T_{1}^{*} \bigoplus T_{2}^{*}$.
 \end{proof}
 \begin{corollary}\label{cor 3.4}
 Let $T_{1}: D(T_{1})\subset H_{1} \to K_{1} \text{ and } T_{2}: D(T_{2}) \subset H_{2} \to K_{2}$ be two densely defined closed operators with closed ranges $R(T_{1})$ and $R(T_{2})$. Then $((T_{1} \bigoplus T_{2})^{\dagger})^{*} = ((T_{1} \bigoplus T_{2})^{*})^{\dagger}$.
 \end{corollary}
 \begin{proof}
 By Theorem \ref{thm 3.1}, Lemma \ref{lemma 3.3} and closed ranges $R(T_{i}^{*}) ~(i= 1,2)$ say that
 \begin{align*}
 ((T_{1} \bigoplus T_{2})^{\dagger})^{*}= (T_{1}^{\dagger} \bigoplus T_{2}^{\dagger})^{*} = (T_{1}^{\dagger})^{*} \bigoplus ( {T_{2}^{\dagger}})^{*} = (T_{1}^{*} \bigoplus T_{2}^{*})^{\dagger} = ((T_{1} \bigoplus T_{2})^{*})^{\dagger}.
 \end{align*}
 \end{proof}
 \begin{remark}
 When $T_{i}: D(T_{i}) \subset H_{i} \to K_{i} ~(i = 1,2, \dots, n)$ are densely defined closed operators with closed ranges $R(T_{i}), i = 1,2, \dots, n$. Then $((T_{1} \bigoplus T_{2} \bigoplus \dots \bigoplus T_{n})^{\dagger})^{*} = ((T_{1} \bigoplus T_{2} \bigoplus \dots \bigoplus T_{n})^{*})^{\dagger}$.
 \end{remark}
 \begin{corollary}
 Let $T_{i}:  D(T_{i}) \subset H_{i} \to K_{i} ~(i= 1, 2)$ be two closed operators with closed ranges $R(T_{i}), ~i = 1,2$. Then $\gamma(T_{1} \bigoplus T_{2}) = \min \{\gamma(T_{1}) ,\gamma(T_{2})\} > 0$, where $\gamma(T)$ is the reduced minimum modulus of $T$.
 \end{corollary}
 \begin{proof}
 Since, $R(T_{1})$ and $R(T_{2})$ are closed. So, $R(T_{1} \bigoplus T_{2})$ is also closed. Thus $\gamma(T_{1} \bigoplus T_{2}) > 0$. 
 Now, for all $(k_{1}, k_{2}) \in K_{1} \bigoplus K_{2}$, we have
 \begin{align*}
 \|(T_{1}^{\dagger} \bigoplus T_{2}^{\dagger})(k_{1}, k_{2})\| &= \|(T_{1}^{\dagger}k_{1}, T_{2}^{\dagger}k_{2})\|\\
 &= (\|T_{1}^{\dagger}k_{1}\|^{2} + \|T_{2}^{\dagger}k_{2}\|^{2})^{\frac{1}{2}}\\
 & \leq (\|T_{1}^{\dagger}\|^{2} \|k_{1}\|^{2} + \|T_{2}^{\dagger}\|^{2}\|K_{2}\|^{2})^{\frac{1}{2}}\\
 &\leq \max\{\|T_{1}^{\dagger}\|, \| T_{2}^{\dagger}\|\} \|(k_{1}, k_{2})\|.
 \end{align*} 
 So,
 \begin{equation}\label{eqn 12}
 \|T_{1}^{\dagger} \bigoplus T_{2}^{\dagger}\| \leq \max\{ \|T_{1}^{\dagger}\|, \|T_{2}^{\dagger}\| \}
 \end{equation}
 Again, for all $w_{1} \in K_{1}$, $\|(T_{1}^{\dagger} \bigoplus T_{2}^{\dagger})(w_{1},0)\| = \|T_{1}^{\dagger}w_{1}\| \leq \|T_{1}^{\dagger} \bigoplus T_{2}^{\dagger}\| \|w_{1}\|$. Similarly,
 For all $w_{2} \in K_{2}$, $\|(T_{1}^{\dagger} \bigoplus T_{2}^{\dagger})(0, w_{2})\| = \|T_{2}^{\dagger}w_{2}\| \leq \|T_{1}^{\dagger} \bigoplus T_{2}^{\dagger}\| \|w_{2}\|$. Hence,
 \begin{equation}\label{eqn 13}
 \max\{ \|T_{1}^{\dagger}\|, \|T_{2}^{\dagger}\|\} \leq \|T_{1}^{\dagger} \bigoplus T_{2}^{\dagger}\|.
 \end{equation}
 The equations (\ref{eqn 12}) and (\ref{eqn 13}) guarantee that $\max\{ \|T_{1}^{\dagger}\|, \|T_{2}^{\dagger}\|\} = \|T_{1}^{\dagger} \bigoplus T_{2}^{\dagger}\|$.
 Therefore, 
 \begin{align*}
 \gamma(T_{1} \bigoplus T_{2}) = \frac{1}{\|T_{1}^{\dagger} \bigoplus T_{2}^{\dagger}\|} = \min\{ \frac{1}{\|T_{1}^{\dagger}\|}, \frac{1}{\|T_{2}^{\dagger}\|}\} = \min\{\gamma(T_{1}), \gamma(T_{2})\} > 0.
 \end{align*}
  \end{proof}
\begin{theorem}\label{thm 3.7}
Let $T_{1}: D(T_{1}) \subset H_{1} \to K_{1}$ and $T_{2}: D(T_{2}) \subset H_{2} \to K_{2}$ be two densely defined closed operators with closed ranges $R(T_{1})$ and $R(T_{2})$. Then $|(T_{1} \bigoplus T_{2})^{\dagger}| = |T_{1}^{\dagger}| \bigoplus |T_{2}^{\dagger}|$.
\end{theorem}
\begin{proof}
By Theorem \ref{thm 3.1}, we have 
\begin{align*}
|(T_{1} \bigoplus T_{2})^{\dagger}| = |T_{1}^{\dagger} \bigoplus T_{2}^{\dagger}| &= ((T_{1}^{\dagger} \bigoplus T_{2}^{\dagger})^{*} (T_{1}^{\dagger} \bigoplus T_{2}^{\dagger}))^{\frac{1}{2}}\\ 
&= ((T_{1}^{*})^{\dagger}T_{1}^{\dagger} \bigoplus (T_{2}^{*})^{\dagger}T_{2}^{\dagger})^{\frac{1}{2}}\\
&= ((T_{1}T_{1}^{*})^{\dagger} \bigoplus (T_{2}T_{2}^{*})^{\dagger} )^{\frac{1}{2}}\\
&= ((|T_{1}^{*}|^{2})^{\dagger} \bigoplus (|T_{2}^{*}|^{2})^{\dagger})^{\frac{1}{2}}
\end{align*}
Now, we claim $(|T_{1}^{*}|^{2})^{\dagger} = (|T_{1}^{*}|^{\dagger})^{2}$.
Since, $R(T_{1})$ is closed. So, $R(|T_{1}^{*}|^{2}) = R(T_{1}T_{1}^{*}) = R(T_{1})$ is closed which implies  $R(|T_{1}^{*}|^{2}) = R(|T_{1}^{*}|)$. Then, $N(|T_{1}^{*}|^{2})^{\perp} = N(|T_{1}^{*}|)^{\perp}$. It is easy to prove that $D(|T_{1}^{*}|^{2} (|T_{1}^{*}|^{\dagger})^{2}) = K_{1} = P_{R(|T_{1}^{*}|^{2})}$ and $D((|T_{1}^{*}|^{\dagger})^{2} |T_{1}^{*}|^{2}) = D(|T_{1}^{*}|^{2})$. So,
\begin{align*}
|T_{1}^{*}|^{2} (|T_{1}^{*}|^{\dagger})^{2} &= |T_{1}^{*}||T_{1}^{*}||T_{1}^{*}|^{\dagger}|T_{1}^{*}|^{\dagger}\\
&= |T_{1}^{*}| P_{R(|T_{1}^{*}|)} |T_{1}^{*}|^{\dagger}\\
&= |T_{1}^{*}||T_{1}^{*}|^{\dagger}\\
&= P_{R(|T_{1}^{*}|)}\\
&= P_{R(|T_{1}^{*}|^{2})}.
\end{align*}
Again,
\begin{align*}
(|T_{1}^{*}|^{\dagger})^{2} |T_{1}^{*}|^{2} &= |T_{1}^{*}|^{\dagger} |T_{1}^{*}|^{\dagger} |T_{1}^{*}| |T_{1}^{*}|\\
&= P_{N(|T_{1}^{*}|)^{\perp}}{\vert_ {D(|T_{1}^{*}|^{2})}}\\
&= P_{N(|T_{1}^{*}|^{2})^{\perp}}{\vert_ {D(|T_{1}^{*}|^{2})}}
\end{align*}
Moreover, $(|T_{1}^{*}|^{\dagger})^{2} |T_{1}^{*}|^{2} (|T_{1}^{*}|^{\dagger})^{2} = |T_{1}^{*}|^{\dagger} P_{N(|T_{1}^{*}|)^{\perp}}{\vert _{D(|T_{1}^{*}|)}} P_{R(|T_{1}^{*}|)}|T_{1}^{*}|^{\dagger} = (|T_{1}^{*}|^{\dagger})^{2}$.
By Theorem 5.7 \cite{MR0451661}, we get $(|T_{1}^{*}|^{\dagger})^{2} = (|T_{1}^{*}|^{2})^{\dagger}$. Similarly, we have $(|T_{2}^{*}|^{\dagger})^{2} = (|T_{2}^{*}|^{2})^{\dagger}$. Hence,
\begin{align*}
|(T_{1} \bigoplus T_{2})^{\dagger}| &= ((|T_{1}^{*}|^{2})^{\dagger} \bigoplus (|T_{2}^{*}|^{2})^{\dagger})^{\frac{1}{2}}\\
&= ((|T_{1}^{*}|^{\dagger})^{2} \bigoplus (|T_{2}^{*}|^{\dagger})^{2})^{\frac{1}{2}}\\
&= (|T_{1}^{*}|^{\dagger}) \bigoplus (|T_{2}^{*}|^{\dagger})\\
&=|T_{1}^{\dagger}| \bigoplus |T_{2}^{\dagger}| ~(\text {by Proposition 3.19 \cite{MR4282727}}).
\end{align*}
\end{proof}
\begin{corollary}
Let $T_{1}: D(T_{1}) \subset H_{1} \to K_{1}$ and $T_{2}: D(T_{2}) \subset H_{2} \to K_{2}$ be two densely defined closed operators with closed ranges $R(T_{1})$ and $R(T_{2})$. Then $|T_{1} \bigoplus T_{2}|^{\dagger} = |T_{1}|^{\dagger} \bigoplus |T_{2}|^{\dagger} = |((T_{1} \bigoplus T_{2})^{*})^{\dagger}|$
\end{corollary}
\begin{proof}
Since $R(T_{1})$ is closed. So, $R(|T_{1}|^{2}) = R(T_{1}^{*}T_{1}) = R(T_{1}^{*})$ which implies $R(|T_{1}|^{2})$ is closed and $R(|T_{1}|^{2}) = R(|T_{1}|)$. Similarly, $R(|T_{2}|^{2}) = R(|T_{2}|)$ is closed.
\begin{align*}
|T_{1} \bigoplus T_{2}|^{\dagger} = (((T_{1} \bigoplus T_{2})^{*} (T_{1} \bigoplus T_{2}))^{\frac{1}{2}})^{\dagger} &= ((T_{1}^{*}T_{1} \bigoplus T_{2}^{*}T_{2})^{\frac{1}{2}})^{\dagger}\\
&= ((|T_{1}|^{2} \bigoplus |T_{2}|^{2})^{\frac{1}{2}})^{\dagger}\\
&= (|T_{1}| \bigoplus |T_{2}|)^{\dagger}\\
&= |T_{1}|^{\dagger} \bigoplus |T_{2}|^{\dagger} ~(\text{by Theorem \ref{thm 3.1}}).
\end{align*}
Now, 
\begin{align*}
|T_{1}|^{\dagger} \bigoplus |T_{2}|^{\dagger} &= |(T_{1}^{*})^{\dagger}| \bigoplus |(T_{2}^{*})^{\dagger}| ~( \text{ by Theorem 3.19 \cite{MR4282727}})\\
&= |(T_{1}^{*} \bigoplus T_{2}^{*})^{\dagger}| ~(\text{ by Theorem \ref{thm 3.7}})\\
&= |((T_{1} \bigoplus T_{2})^{*})^{\dagger}| ~(\text{ by Lemma \ref{lemma 3.3}})
\end{align*}
Therefore, $|T_{1} \bigoplus T_{2}|^{\dagger} = |T_{1}|^{\dagger} \bigoplus |T_{2}|^{\dagger} = |((T_{1} \bigoplus T_{2})^{*})^{\dagger}|$.
\end{proof}
\begin{example}
Let us define $T_{1}: \ell^{2} \to \ell^{2}$ by $T_{1}(x_{1}, x_{2}, \dots, x_{n}, \dots) = (x_{1}, 2x_{2}, \dots, nx_{n}, \dots)$ and $T_{2}:\ell^{2} \to \ell^{2}$ by $T_{2}(x_{1}, x_{2}, \dots, x_{n}, \dots) = (0, 2x_{2}, \dots, nx_{n}, \dots)$. Then $T_{1} = T_{1}^{*}$ with $R(T_{1}) = \ell^{2}$ is closed (Example 2.4 \cite{majumdar2024hyers}). Again, $T_{2} = T_{2}^{*}$ with $R(T_{2}) = \{(0, y_{2}, y_{3}, \dots, y_{n}, \dots ) : (0, y_{2}, y_{3}, \dots, y_{n}, \dots ) \in \ell^{2}\}$ is closed (Example 2.19 \cite{majumdar2024hyers}). Theorem \ref{thm 3.1} says that $(T_{1} \bigoplus T_{2})^{\dagger} = T_{1}^{\dagger} \bigoplus T_{2}^{\dagger}$. Here,
\begin{align*}
T_{1}^{\dagger} (x_{1}, x_{2}, \dots, x_{n}, \dots) = (x_{1}, \frac{x_{2}}{2}, \dots, \frac{x_{n}}{n}, \dots) ,\text{ for all } (x_{1}, x_{2},\dots, x_{n}, \dots) \in \ell^{2}.
\end{align*}
and
\begin{align*}
& T_{2}^{\dagger}(0, x_{2}, \dots, x_{n}, \dots) = (0, \frac{x_{2}}{2}, \dots, \frac{x_{n}}{n}, \dots), \text{ for all}~ (0, x_{2}, \dots, x_{n}, \dots) \in R(T_{2})\\
& T_{2}^{\dagger}(x_{1},0, \dots, 0, \dots) = (0, 0, \dots, 0, \dots), \text{ for all } (x_{1},0, \dots,0, \dots)\in R(T_{2})^{\perp}
\end{align*}
Then, $\text{ for all } (w,z) = ((w_{1}, w_{2},\dots, w_{n},\dots), (0,z_{1},\dots,z_{n},\dots))\in R(T_{1}) \bigoplus R(T_{2})$, we have 
\begin{align*}
(T_{1} \bigoplus T_{2})^{\dagger} (w,z) &=(T_{1} \bigoplus T_{2})^{\dagger}((w_{1}, w_{2},\dots, w_{n},\dots), (0,z_{1},\dots,z_{n},\dots))\\
&= ((w_{1},\frac{w_{2}}{2}, \dots, \frac{w_{n}}{n},\dots), (0,\frac{z_{2}}{2}, \dots, \frac{z_{n}}{n},\dots))
\end{align*}
with
\begin{align*}
    (T_{1} \bigoplus T_{2})^{\dagger} (u,v) = (0,0), \text{ for all } (u,v) \in (R(T_{1}) \bigoplus R(T_{2}))^{\perp}.
\end{align*}
\end{example}
In the next theorem, we present the Moore-Penrose inverse of the sum of two operators.
\begin{theorem}
Let $S \in B(H, K)$ and $T$ be a densely defined closed operator from $H$ into $K$ with the following conditions:
\begin{enumerate}
\item $\|T^{\dagger}S\| < 1$.
\item $\|Sx\| \leq b\|Tx\| , \text{ for all } x\in D(T) ~\text{and}~ 0<b<1$.
\item $\|S^{*}x\| \leq c \|T^{*}x\|, \text{ for all } x\in D(T^{*}) \text{ and } 0 <c<1$.
\end{enumerate}
Then $(T + S)^{\dagger} = (I + T^{\dagger}S)^{-1}T^{\dagger}$.
\end{theorem}
\begin{proof}
From the given condition (3), we get $\langle SS^{*}x, x \rangle \leq c  \langle TT^{*}x, x \rangle , \text{ for all } x \in D(TT^{*})$. By Theorem 2 \cite{MR0203464}  and $ 0 < c < 1$, we get a contraction operator $U$ such that $S = TU$ (Since, $D(S) = H$). Then $R(S) \subset R(T)$ which justifies the boundedness of $T^{\dagger}S$. By condition (1), we have the existence of the bounded operator $(I + T^{\dagger}S)^{-1}$ in $B(H)$.
We will show that the domain of $(T+S)$ is decomposable.

Condition (2) says that $N(T + S) = N(T)$. Again, $C(T + S) = N(T + S)^{\perp} \cap D(T + S) = N(T)^{\perp} \cap D(T) = C(T)$. So, 
\begin{align*}
D(T + S) = D(T) = N(T) + C(T) = N(T + S) + C(T + S).
\end{align*}
Thus, $D(T + S)$ is decomposable and the Moore-Penrose inverse of $(T + S)$ exists.
Let us consider $v \in R(T + S)$, then there exists $v_{1} \in H$ such that $v = (T + S)v_{1} = (T + TU)v_{1} \in R(T)$. Thus, $R(T + S) \subset R(T)$. Considering an arbitrary element $p \in R(T)$, then there are two elements $p_{1}, p_{2} \in H$ with $p_{1} = (I + T^{\dagger}S)p_{2}$ such that 
\begin{align*}
p = Tp_{1} = T(I + T^{\dagger}S)p_{2} = (T + S)p_{2} \in R(T + S).
\end{align*}
Hence, $R(T) = R(T + S)$. It is true that $TT^{\dagger}S = S ~( \text{because } R(S) \subset R(T))$. Moreover, $z = z_{1} + z_{2}\in D(T^{\dagger}T)$ where $z_{1} \in N(T) \text{ and } z_{2} \in C(T)$, we have
\begin{align*}
ST^{\dagger}Tz = S{P_{N(T)^{\perp}}}{\vert_{D(T)}z }= Sz_{2} = Sz
\end{align*}
which shows that $ ST^{\dagger}T = S$ on domain $D(T^{\dagger}T)$. Now,
\begin{align*}
(T + S)(I + T^{\dagger}S)^{-1}T^{\dagger} &= (T + TT^{\dagger}S)(I + T^{\dagger}S)^{-1}T^{\dagger}\\
&= TT^{\dagger}\\
&= P_{R(T)}{\vert_{ D(T^{\dagger})}}\\
&= P_{R(T + S)}{\vert_{(D(T + S)^{\dagger})}}.
\end{align*}
Similarly,
\begin{align*}
(I + T^{\dagger}S)^{-1}T^{\dagger}(T + S) &= (I + T^{\dagger}S)^{-1} (T^{\dagger}T + T^{\dagger}ST^{\dagger}T)\\
&= T^{\dagger}T\\
&= P_{N(T)^{\perp}}{\vert_{D(T)}}\\ 
&= P_{N(T + S)^{\perp}}{\vert_{D(T + S)}}.
\end{align*}
Therefore, the Moore-Penrose inverse of $(T + S)$ is $(I + T^{\dagger}S)^{-1}T^{\dagger}$.

\end{proof}

\begin{center}
	\textbf{Acknowledgements}
\end{center}

\noindent The present work of the second author was partially supported by Science and Engineering Research Board (SERB), Department of Science and Technology, Government of India (Reference Number: MTR/2023/000471) under the scheme ``Mathematical Research Impact Centric Support (MATRICS)''. 

\end{document}